\documentclass[12pt,twoside]{article}
\usepackage{amssymb,amsmath} % font used for R in Real numbers
\newtheorem{theorem}{Theorem}[section]
\newtheorem{corollary}[theorem]{Corollary}
\newtheorem{lemma}[theorem]{Lemma}

\newtheorem{conjecutre}[theorem]{Conjecutre}
\numberwithin{equation}{section}
\newcommand \proof{\noindent\textbf{Proof. \ }}
\newcommand \epf{\hfill \rule{0.2cm}{0.2cm}\bigskip \par}

\title{Asymptotic behavior of extremal functions to an inequality involving Hardy potential and critical Sobolev exponent}

\begin{document}
\author{{Benjin Xuan}\thanks{Supported by Grant 10101024 and 10371116 from
the National Natural Science Foundation of China. {\it e-mail:bjxuan@ustc.edu.cn}\,(B. Xuan)}\ \ \  \ Jiangchao Wang\\
{\it Department of Mathematics}\\
{\it University of Science and Technology of China}\\
{\it Universidad Nacional de Colombia}\\
}
\date{}
\maketitle

\begin{abstract}
In this paper, we study the asymptotic behavior of radial extremal
functions to an inequality involving Hardy potential and
critical Sobolev exponent. Based on the asymptotic behavior at the origin and the infinity,
we shall deduce a strict inequality
between two best constants. Finally, as an application of this strict inequality, we consider the existence of
nontrivial solution of a quasilinear Brezis-Nirenberg type problem with Hardy potential and critical Sobolev exponent.


\end{abstract}

\section{Introduction.}\label{intro}

In this paper, we study the asymptotic behavior of extremal
functions to the following inequality involving Hardy potential and
critical Sobolev exponent:
\begin{equation}\label{eq1.1}
C \left(\int_{\mathbb{R}^N} \frac{|u|^{p_*}}{|x|^{bp_*}}
\,{\rm d}x \right)^{p/p_*} \leqslant  \int_{\mathbb{R}^N} \left(
\frac{|{\rm D}u|^{p}}{|x|^{ap}} -\mu
\frac{|u|^{p}}{|x|^{(a+1)p}}\right)\,{\rm d}x,
\end{equation}
where $1<p<N,\ 0\leqslant a <\frac{N-p}{p},\ a\leqslant b<(a+1),\
p_*=\frac{Np}{N-(a+1-b)p},\ \mu< \overline{\mu}$,
$\overline{\mu}$ is the best constant in the Hardy equality.

We shall show that for $\mu<\overline{\mu}$ the best constant of
inequality (\ref{eq1.1}) is achievable. Furthermore, the extremal functions of inequality (\ref{eq1.1}) is radial
symmetric. Then we study the asymptotic behavior of the radial extremal
functions of inequality (\ref{eq1.1}) at the origin and the
infinity. At last, for any smooth bounded open domain $\Omega
\subset \mathbb{R}^N$ containing $0$ in its insides, we shall
deduce a strict inequality between two best constants
$S_{\lambda,\,\mu}(p,a,b;\Omega)$ and
$S_{0,\,\mu}(p,a,b;\Omega)=S_{0,\,\mu}$:
\begin{equation}
\label{eq1.4}  S_{\lambda,\,\mu}(p,a,b;\Omega)< S_{0,\,\mu},
\end{equation}
if $\lambda>0$, where $S_{0,\,\mu}$ and $S_{\lambda,\,\mu}(p,a,b;\Omega)$ will be
defined in Section 2 and 4 respectively. We believe that the
strict inequality (\ref{eq1.4}) will be useful to study the
existence of quasilinear elliptic problem involving Hardy
potential and critical Sobolev exponent.
As an application of this strict inequality, we consider the existence of
nontrivial solution of a quasilinear Brezis-Nirenberg problem
with Hardy potential and critical Sobolev exponent.

In their famous paper \cite{BN}, Brezis and Nirenberg
studied problem:
\begin{equation}\label{eq1.5}
\left\{ \begin{array}{rl}
-\Delta u &=\lambda u+u^{2^\ast-1}, \mbox{  in }\ \Omega,\\
u& >0,\mbox{  in }\ \Omega,\\
u& =0,  \mbox{  on }\
\partial\,\Omega.
\end{array}\right.
\end{equation}
Since the embedding $H_0^1(\Omega) \hookrightarrow
L^{2^*}(\Omega)$ is not compact where $2^*=2N/(N-2)$, the
asociated energy functional does not satisfy the (PS)
condition globally, which caused a serious difficulty when trying
to apply standard variational methods. Brezis and Nirenberg
successfully reduced the existence of solutions of problem
(\ref{eq1.5}) into the verification of a special version of the
strict inequality (\ref{eq1.4}) with $p=2,\ a=b=\mu=0$. To verify
(\ref{eq1.4}) in their case, they applied the explicit expression
of the extremal functions to the Sobolev inequality, especially
the asymptotic behavior of the extremal functions at the origin
and the infinity. Brezis-Nirenberg type problems have been
generalized to many other situations (see \cite{CG, EH1, EH2, GV, JS,
NL, PS, XC, XSY,ZXP} and references therein).

Recently, Jannelli \cite{JE} introduced the term $\mu\frac{u}{|x|^2}$ in the equation, that is,
\begin{equation}\label{eq1.6}
\left\{ \begin{array}{rl}
-\Delta u - \mu\frac{u}{|x|^2} & =\lambda u+u^{2^\ast-1}, \mbox{  in }\ \Omega,\\
u& >0,\mbox{  in }\ \Omega,\\
u& =0,  \mbox{  on }\
\partial\,\Omega.
\end{array}\right.
\end{equation}
He studied the relation between critical dimensions for $\lambda\in
(0,\ \lambda_1)$ and $L_{\rm loc}^2$ integrability of the
associated Green function, where $\lambda_1$ is the first
eigenvalue of operator $-\Delta - \mu\frac{1}{|x|^2}$ on $\Omega$
with zero-Dirichlet condition. Ruiz and Willem \cite{RW} also
studied problem (\ref{eq1.6}) under various assumption on the
domain $\Omega$, and even for $\mu\leqslant 0$. Those proofs in
\cite{JE} and \cite{RW} were reduced to verify the strict
inequality (\ref{eq1.4}) with $p=2, a=b=0$. In 2001, Ferrero and
Gazzola \cite{FG} considered the existence of sign-changed
solution to problem (\ref{eq1.6}) for larger $\lambda$. They
distinguished two distinct cases: resonant case and non-resonant
cases of the Brezis-Nirenberg type problem (\ref{eq1.6}). For the
resonant case, they only studied a special case: $\Omega$ is the
unit ball and $\lambda=\lambda_1$. The general case was left as an
open problem. In 2004, Cao and Han \cite{CH} complished the
general case. In all the references cited above, the asymptotic
behavior of the extremal functions at the origin and the infinity
was applied to derived the local (PS) condition for the associated
energy functional.

The rest of this paper is organized as follows. In section 2, we shall show that the best constant of (\ref{eq1.1})
is achieved by some radial extremal functions.
Section 3 is concerning with the asymptotic behavior of the radial extremal functions. In Section 4, we first derive various estimates on
the approximation extremal functions, and then establish the strict inequality (\ref{eq1.4}).
In section 5, based on this strict inequality, we obtain the existence  results of
nontrivial solution of a quasilinear Brezis-Nirenberg problem.

%%%%%%%%%%%%%%%%%%%%%%%%%%%%%%%%%%%%%%%%%%%%%%%%%%%%%%%%%%%%%%%%%%%%%%%
\section{Radial extremal functions}\label{radial}
%%%%%%%%%%%%%%%%%%%%%%%%%%%%%%%%%%%%%%%%%%%%%%%%%%%%%%%%%%%%%%%%%%%%%%%%%%

In order to obtain the extremal functions of (\ref{eq1.1}). We consider the following extremal problem:
\begin{equation}\label{eq2.2}
S_{0,\,\mu}=\inf\left\{Q_\mu(u)\ :\ u\in
\mathfrak{D}_{a,b}^{1,p}(\mathbb{R}^N),\
 \|u;L_b^{p_\ast}(\mathbb{R}^N)\|=1\right\},
\end{equation}
where
$$
Q_\mu(u)=\int_{\mathbb{R}^N}\frac{|{\rm D}u|^p}{|x|^{ap}}\,{\rm
d}x-\mu\int_{\mathbb{R}^N}\frac{|u|^p}{|x|^{(a+1)p}}\,{\rm d}x,
$$
and
$$
\mathfrak{D}_{a,b}^{1,p}(\mathbb{R}^N)=\{u\in
L_b^{p_\ast}(\mathbb{R}^N):|{\rm D}u|\in L_a^p(\mathbb{R}^N)\}
$$
is the closure of $C_0^\infty (\mathbb{R}^N)$ under the norm $\| u
\|_{\mathfrak{D}_{a,b}^{1,p}(\mathbb{R}^N)}=\| |{\rm D}u|; L_a^p
(\mathbb{R}^N) \|$. For any $\alpha, q$, the norm of weighted space $L_\alpha^{q}(\mathbb{R}^N)$ is defined as
$$
\|u;L_\alpha^{q}(\mathbb{R}^N)\|=(\int_{\mathbb{R}^N}\frac{|u|^{q}}{|x|^{\alpha{q}}}\,{\rm
d}x)^{\frac{1}{q}}.
$$

Similar to Lemma 2.1 in \cite{GP}, one can easily obtain the
following Hardy inequality with best constant
$\overline{\mu}=(\frac{N-(a+1)p}{p})^p$:
\begin{equation}\label{eq2.02}
\overline{\mu} \int_{\mathbb{R}^N}\frac{|u|^p}{|x|^{(a+1)p}}\,{\rm d}x \leqslant \int_{\mathbb{R}^N}\frac{|{\rm D}u|^p}{|x|^{ap}}\,{\rm
d}x.
\end{equation}
Thus, for $\mu<\overline{\mu}$, $Q_\mu(u)\geqslant 0$ for all
$u\in\mathfrak{D}_{a,b}^{1,p}(\mathbb{R}^N)$, and the equality holds if
and only if $u\equiv 0$. From the so-called
Caffarelli-Kohn-Nirenberg inequality \cite{CKN},
$S_{0,\,\mu}<\infty$.

\begin{lemma}\label{lem2.3} If
$\mu\in (0,\ \overline{\mu}),\ b\in [a,\ a+1)$, then $S_{0,\,\mu}$ is
achieved at some nonnegative function $u_0\in
\mathfrak{D}_{a,b}^{1,p}(\mathbb{R}^N)$. In particular, there
exists a solution to the following ``limited equation":
\begin{equation}\label{eq2.1}
-\mbox{div}(\frac{{|{\rm D}u|}^{p-2}{\rm
D}u}{|x|^{ap}})-\mu\frac{|u|^{p-2}u}{|x|^{(a+1)p}}=\frac{|u|^{p_\ast-2}u}{|x|^{bp_\ast}}.
\end{equation}
 \end{lemma}

\proof  The achievability of $S_{0,\,\mu}$ at some $u_0\in
\mathfrak{D}_{a,b}^{1,p}(\mathbb{R}^N)$ with $\|u_0;L_b^{p_\ast}(\mathbb{R}^N)\|=1$ is due to \cite{WW} for
$p=2$ and \cite{TY} for general $p$. Without loss of generality,
suppose that $u_0\geqslant 0$, otherwise, replace it by $|u_0|$. It
is easy to see that $u_0$ satisfies the following Euler-Lagrange
equation:
$$
-\mbox{div}(\frac{|{\rm D}u|^{p-2}{\rm
D}u}{|x|^{ap}})-\mu\frac{|u|^{p-2}u}{|x|^{(a+1)p}}
=\delta\frac{|u|^{p_\ast-2}u}{|x|^{b{p_\ast}}},
$$
where $\delta=Q_\mu(u_0)/ \|u_0;L_b^{p_\ast}(\mathbb{R}^N)\|^{p_\ast}= Q_\mu(u_0)=S_{0,\,\mu}>0$ is the Lagrange multiplier. Set
$\overline{u}=c_0u_0,\ c_0=S_{0,\,\mu}^{\frac{1}{p_\ast-p}}$, then
$\overline{u}$ is a solution to equation (\ref{eq2.1}). \epf

In fact, all the dilation of $u_0$ of the form
$\sigma^{-\frac{N-(a+1)p}{p}}u_0(\frac{\cdot}{\sigma})$ are also
minimizers of $S_{0,\,\mu}$. In order to obtain further
properties of the minimizers of $S_{0,\,\mu}$, let's recall the
definition of the Schwarz symmetrization (see \cite{HT}). Suppose
that $\Omega\subset{\mathbb{R}^N}$, and $f\in C_0(\Omega)$ is a
nonnegative continuous function with compact support, the the
Schwarz symmetrization $S(f)$ of $f$ is defined as
$$
S(f)(x)=\sup \left\{t:\mu(t)>\omega_N|x|^N\right\},
\quad \mu(t)=\left|\,\{x\ :\ f(x)>t\}\right|,$$
 where $\omega_N$ denotes the volume of the standard $N$-sphere.
Applying those properties of Schwarz symmetrization in \cite{HT}, we have the following lemma:
\begin{lemma}\label{lem2.2} For
$v\in\mathfrak{D}_{a,b}^{1,p}(\mathbb{R}^N)\setminus \{0\},\
k\geqslant 0$, define
$$
R(v)=\dfrac{\int_{S^{N-1}}\int_0^{+\infty}\{k^{1-p-\frac{p}{p_\ast}}(|\partial_\rho
v|^2+\frac{|\Lambda v|^2}{\rho^2})^{p/2}\rho^{N-1}
-k^{1-\frac{p}{p_\ast}}\mu|v|^p\rho^{N-1-p}\}\,{\rm d}\rho \,{\rm
d}S}{\int_{S^{N-1}}\int_0^{+\infty}|v|^{p_\ast}\rho^{\frac{(N-p){p_\ast}}{p}-1}\,{\rm
d}\rho \,{\rm d}S},
$$
where $\partial_\rho$ is the directional differential operator
along direction $\rho$ and $\Lambda$ is the tangential
differential operator on $S^{N-1}$. Then
$$
\inf \{R(v)\ :\ v\in\mathfrak{D}_{a,b}^{1,p}(\mathbb{R}^N)  \mbox{
is radial}\} =  \inf \{R(v)\ :\ v
\in\mathfrak{D}_{a,b}^{1,p}(\mathbb{R}^N)\}.
$$
 \end{lemma}
\proof By the density argument, it suffices to prove the lemma for
$v\in C_0^\infty (\mathbb{R}^N)$. Let $v^*$ be the Schwarz
symmetrization of $v$. Noting that $\Lambda v^\ast=0, p_*\leqslant \frac{Np}{N-p}$, and applying those properties of Schwarz symmetrization in \cite{HT}, we have
$$
\int_{S^{N-1}}\int_0^{+\infty}|v^\ast|^{p_\ast}\rho^{\frac{(N-p){p_\ast}}{p}-1} \,{\rm d}\rho
 \,{\rm d}S\geqslant \int_{S^{N-1}}\int_0^{+\infty}|v|^{p_\ast}\rho^{\frac{(N-p){p_\ast}}{p}-1} \,{\rm d}\rho
 \,{\rm d}S=1,
$$
\begin{equation*}
\begin{split}
k^{1-p-\frac{p}{p_\ast}}\int_{S^{N-1}}\int_0^{+\infty}&(|\partial_\rho
v^\ast|^2+\frac{|\Lambda v^\ast|^2}{\rho^2})^{p/2}\rho^{N-1}\,{\rm
d}\rho
\,{\rm d}S\\
&\leqslant
k^{1-p-\frac{p}{p_\ast}}\int_{S^{N-1}}\int_0^{+\infty}(|\partial_\rho
v|^2+\frac{|\Lambda v|^2}{\rho^2})^{p/2}\rho^{N-1}\,{\rm d}\rho
\,{\rm d}S
\end{split}
\end{equation*}
and
$$k^{1-\frac{p}{p_\ast}}\mu\int_{S^{N-1}}\int_0^{+\infty}|v^\ast|^p\rho^{N-1-p}\,{\rm d}\rho
\,{\rm d}S\geqslant
k^{1-\frac{p}{p_\ast}}\mu\int_{S^{N-1}}\int_0^{+\infty}|v|^p\rho^{N-1-p}\,{\rm
d}\rho \,{\rm d}S.
$$
Thus, we have
\begin{equation*}
\begin{split}
&\int_{S^{N-1}}\int_0^{+\infty}\{k^{1-p-\frac{p}{p_\ast}}(|\partial_\rho
v^\ast|^2+\frac{|\Lambda v^\ast|^2}{\rho^2})^{p/2}\rho^{N-1}
-k^{1-\frac{p}{p_\ast}}\mu|v^\ast|^p\rho^{N-1-p}\}\,{\rm d}\rho
\,{\rm d}S\\&\leqslant
\int_{S^{N-1}}\int_0^{+\infty}\{k^{1-p-\frac{p}{p_\ast}}(|\partial_\rho
v|^2+\frac{|\Lambda v|^2}{\rho^2})^{p/2}\rho^{N-1}
-k^{1-\frac{p}{p_\ast}}\mu |v|^p\rho^{N-1-p}\}\,{\rm d}\rho \,{\rm
d}S.
\end{split}
\end{equation*}
That is,
$$
R(v^*)\leqslant R(v),
$$
thus,
$$
\inf \{R(v)\ :\ v\in\mathfrak{D}_{a,b}^{1,p}(\mathbb{R}^N)  \mbox{
is radial}\} \leqslant  \inf \{R(v)\ :\ v
\in\mathfrak{D}_{a,b}^{1,p}(\mathbb{R}^N)\}.
$$
On the other hand, it is trivial that
$$
\inf \{R(v)\ :\ v\in\mathfrak{D}_{a,b}^{1,p}(\mathbb{R}^N)  \mbox{
is radial}\} \geqslant  \inf \{R(v)\ :\ v
\in\mathfrak{D}_{a,b}^{1,p}(\mathbb{R}^N)\}.
$$
\epf

\begin{lemma}\label{lem2.4} If
$\mu\in (0,\ \overline{\mu}),\ b\in [a,\ a+1)$, then all the minimizers
of $S_{0,\,\mu}$ is radial. In particular, there exists a family of
radial solutions to equation (\ref{eq2.1}).
 \end{lemma}

\proof We rewrite those integrals in $S_{0,\,\mu}$ in polar
coordinates. Noting that $|{\rm D}u|^2=|\partial_r
u|^2+\frac{1}{r^2}|\Lambda u|^2$, we have
\begin{equation}\label{eq2.6}
\begin{split}
Q_\mu(u)&=\int_{S^{N-1}}\int_0^{+\infty}(|\partial_r
u|^2+\frac{1}{r^2}|\Lambda u|^2)^{p/2} r^{N-1-ap}\,{\rm d}r\,{\rm d}S\\
&\ \ \ \  -\mu\int_{S^{N-1}}\int_0^{+\infty}|u|^p
r^{N-1-(a+1)p}\,{\rm d}r\,{\rm d}S.
\end{split}
\end{equation}
Making the change of variables $r=\rho^k,\quad
k=\frac{N-p}{N-(a+1)p}\geqslant 1$, from (\ref{eq2.6}), we have
\begin{equation}\label{eq2.7}
\begin{split}
Q_\mu(u)&=k^{1-p}\int_{S^{N-1}}\int_0^{+\infty}(|\partial_\rho
u|^2+k^2\frac{|\Lambda u|^2}{\rho^2})^{p/2}\rho^{N-1}\,{\rm d}\rho
\,{\rm d}S \\
&\ \ \ \  -k\mu\int_{S^{N-1}}\int_0^{+\infty}|u|^p
\rho^{N-p-1}\,{\rm d}\rho \,{\rm d}S.
\end{split}
\end{equation}
On the other hand, the restriction condition
$\|u;L_b^{p_\ast}(\mathbb{R}^N)\|=1$ becomes
\begin{equation}\label{eq2.8}
k\int_{S^{N-1}}\int_0^{+\infty}|u|^{p_\ast}\rho^{\frac{(N-p){p_\ast}}{p}-1}\,{\rm
d}\rho \,{\rm d}S=1.
\end{equation}
To cancel the coefficient $k$ in (\ref{eq2.8}), let
$v=k^{\frac{1}{p_\ast}}u$, then we have the following equivalent
form of $S_{0,\,\mu}$:
\begin{equation}\label{eq2.9}
\begin{split}
 S_{0,\,\mu}= \inf\left
\{ \right. & \left. \int_{S^{N-1}}\int_0^{+\infty}\{k^{1-p-\frac{p}{p_\ast}}(|\partial_\rho
v|^2+k^2\frac{|\Lambda v|^2}{\rho^2})^{p/2}\rho^{N-1} \right.\\
& \ \ \  -k^{1-\frac{p}{p_\ast}}\mu|v|^p\rho^{N-1-p}\}\,{\rm
d}\rho \,{\rm d}S  \ :\ v \in\mathfrak{D}_{a,b}^{1,p}(\mathbb{R}^N),\\
& \left. \ \ \  \ \ \  \ \ \  \ \ \
\int_{S^{N-1}}\int_0^{+\infty}|v|^{p_\ast}\rho^{\frac{(N-p){p_\ast}}{p}-1}\,{\rm
d}\rho \,{\rm d}S=1\right\}.
\end{split}
\end{equation}
Since $k\geqslant 1$, we have
\begin{equation}\label{eq2.10}
\begin{split}
S_{0,\,\mu} \geqslant \inf\left
\{ \right. & \left. \int_{S^{N-1}}\int_0^{+\infty}\{k^{1-p-\frac{p}{p_\ast}}(|\partial_\rho
v|^2+\frac{|\Lambda v|^2}{\rho^2})^{p/2}\rho^{N-1} \right.\\
& \ \ \  -k^{1-\frac{p}{p_\ast}}\mu|v|^p\rho^{N-1-p}\}\,{\rm
d}\rho \,{\rm d}S\  :\  v
\in\mathfrak{D}_{a,b}^{1,p}(\mathbb{R}^N),\\
& \left. \ \ \  \ \ \  \ \ \  \ \ \
\int_{S^{N-1}}\int_0^{+\infty}|v|^{p_\ast}\rho^{\frac{(N-p){p_\ast}}{p}-1}\,{\rm
d}\rho \,{\rm d}S=1\right\}.
\end{split}
\end{equation}
From Lemma \ref{lem2.2}, we know that the left side hand is
achieved at some radial function, and the inequality in
(\ref{eq2.10}) becomes equality if and only if $v$ is radial.
Thus, all the minimizers of $S_{0,\,\mu}$ is radial. \epf

%%%%%%%%%%%%%%%%%%%%%%%%%%%%%%%%%%%%%%%%%%%%%%%%%%%%%%%%%%%%%%%%%%%%%%%
\section{Asymptotic behavior of extremal functions}\label{Behavior}
%%%%%%%%%%%%%%%%%%%%%%%%%%%%%%%%%%%%%%%%%%%%%%%%%%%%%%%%%%%%%%%%%%%%%%%%%%

In this section, we describe the asymptotic behavior of radial
extremal functions of $S_{0,\,\mu}$. Our argument here is similar to
that in \S3.2 of \cite{AFP}.  Let $u(r)$ be a nonnegative radial
solution to (\ref{eq2.1}). Rewriting in polar coordinates, we have
\begin{equation}\label{eq3.1}
(r^{N-1-ap}|u'|^{p-2}u')'+r^{N-1}(\mu\frac{|u|^{p-2}u}{r^{(a+1)p}}+\frac{|u|^{p_\ast-2}u}{r^{bp_\ast}})=0.
\end{equation}
Set
\begin{equation}\label{eq3.2}
t=\log r,\quad y(t)=r^\delta u(r),\quad
z(t)=r^{(1+\delta)(p-1)}|u'(r)|^{p-2}u'(r),
\end{equation}
where $\delta=\frac{N-(a+1)p}{p}$. A simple calculation shows that
\begin{equation}\label{eq3.3}
\left\{{\begin{array}{l} \dfrac{{\rm d}y}{{\rm d}t}=\delta
y+|z|^{\frac{2-p}{p-1}}z;\\[3mm]
\dfrac{{\rm d}z}{{\rm d}t}=-\delta z-|y|^{p_\ast-2}y-\mu|y|^{p-2}y.
\end{array}}\right.
\end{equation}
It follows from (\ref{eq3.3}) that $y$ satisfies the following
equation:
\begin{equation}\label{eq3.4}
(p-1)|\delta y-y'|^{p-2}(\delta y'-y'' )+\delta |\delta
y-y'|^{p-2}(\delta y-y')-\mu y^{p-1}-y^{p_\ast-1}=0.
\end{equation}

 It is easy to see that the complete integral of the
autonomous system (\ref{eq3.3}) is
\begin{equation}\label{eq3.5}
V(y,z)=\frac{1}{p_\ast}|y|^{p_\ast}+\frac{\mu}{p}|y|^p+\frac{p-1}{p}|z|^{\frac{p}{p-1}}+\delta
yz.
\end{equation}
Similar to Lemma 3.6-3.9 in \cite{AFP}, we have the following four
lemmas. We will omit proofs of the first three lemmas because one
only needs to replace $\delta=\frac{N-p}p$ there by
$\delta=\frac{N-(a+1)p}{p}$ in our case. The interested reader can
refer to \cite{AFP}. The idea of the fourth Lemma is also similar
to that of Lemma 3.9 in \cite{AFP}, with different choice of
function $\xi$. We shall write down its complete proof for
completeness.

\begin{lemma}\label{lem3.1} $y$ and $z$ are bounded.
 \end{lemma}
\begin{lemma}\label{lem3.2} For any $t\in\mathbb{R}^N$,
$(y(t),z(t))\in \{(y,z)\in\mathbb{R}^2:V(y,z)=0\}.$
 \end{lemma}
\begin{lemma}\label{lem3.3} There exists $t_0\in\mathbb{R}$, such that $y(t)$ is strictly increasing for $t<t_0$;
and strictly decreasing for $t>t_0$. Furthermore, we have
\begin{equation}\label{eq3.6}
\underset{t\in\mathbb{R}}{\max}\
y(t)=y(t_0)=[\frac{N}{N-(a+1-b)p}(\delta^p-\mu)]^{\frac{1}{p_\ast-p}}
\end{equation}
 \end{lemma}

\begin{lemma}\label{lem3.4} Suppose that $y$ is a positive solution to (\ref{eq3.4}) such that $y$ is increasing
in $(-\infty,0)$ and decreasing in $(0,+\infty)$, then there exist
$c_1,c_2>0$, such that
\begin{equation}\label{eq3.7}
\underset{t\rightarrow-\infty}{\lim}e^{(l_1-\delta)t}y(t)=y(0)c_1>0;
\end{equation}
\begin{equation}\label{eq3.8}
\underset{t\rightarrow+\infty}{\lim}e^{(l_2-\delta)t}y(t)=y(0)c_2>0,
\end{equation}
where $l_1,l_2$ are zeros of function
$\xi(s)=(p-1)s^p-(N-(a+1)p)s^{p-1}+\mu$ such that $0<l_1<l_2.$
 \end{lemma}

\proof First, it is easy to see that $l_1<\delta<l_2$. Next, we
prove (\ref{eq3.7}) step by step and omit the proof of
(\ref{eq3.8}).

\noindent\textbf{1.} It follows from (\ref{eq3.3}) that
\begin{equation}\label{eq3.9}
\begin{array}{rl}
\frac{{\rm d}}{{\rm d}t}(e^{-(\delta-l_1)t}y(t))
&=-(\delta-l_1)e^{-(\delta-l_1)t}y(t)+e^{-(\delta-l_1)t}(\delta y(t)+|z|^{\frac{1}{p-1}})\\
&=e^{-(\delta-l_1)t}y(t)(l_1-\frac{|z(t)|^{\frac{1}{p-1}}}{y(t)}).\\
\end{array}
\end{equation}
Rewritting the above equation into the integral form, we have
\begin{equation}\label{eq3.10}
e^{-(\delta-l_1)t}y(t)=y(0)e^{-\int_{t}^{0}(l_1-y(s)^{-1}|z(s)|^{1/p-1}){\rm
d}s}.
\end{equation}

\noindent\textbf{2.} Let
$H(s)=\frac{|z(s)|^{\frac{1}{p-1}}}{y(s)}$.

\textbf{Claim: }
$H(s)$ is a increasing function from $(-\infty,0]$ into
$(l_1,\delta]$.

In fact, we shall prove that $H'(s)>0$ for $s<0$. Otherwise, we
prove by contradiction, suppose that  there exists $s_0<0$ such
that $H'(s_0)\leqslant0$. A direct computation shows that
$$
H'(s)=\frac{-\frac{1}{p-1}y(s)z'(s)|z(s)|^{\frac{2-p}{p-1}}-|z(s)|^{\frac{1}{p-1}}y'(s)}{y^2(s)}.
$$
Replacing formulas of $y'(s_0)$ and $z'(s_0)$ from (\ref{eq3.3}),
and noting that (\ref{eq3.5}) and Lemma \ref{lem3.2}, it follows
that
$$
H'(s_0)=(\frac{1}{p}-\frac{1}{p_\ast})y^{p_\ast}(s_0)\leqslant0,
$$
which contradicts to the fact that $y>0$. Thus, $H'(s)>0$, and
hence $H$ is strictly increasing on $(-\infty,0]$.

On the other hand, from (\ref{eq3.3}) and $y'(0)=0$, we have
$H(0)=\delta$; from (\ref{eq3.5}), it follows that
$\underset{t\rightarrow-\infty}{\lim}H(s)=l_1$, which proves our
claim.

\noindent\textbf{3.} (\ref{eq3.7}) holds.

From the above claim and (\ref{eq3.10}), it follows that
$e^{-(\delta-l_1)t}y(t)>0$ is decreasing on $(-\infty,0]$, and
hence the limit
$\underset{t\rightarrow-\infty}{\lim}e^{-(\delta-l_1)t}y(t)$
exists. Set
$$
\alpha\equiv\underset{t\rightarrow-\infty}{\lim}e^{-(\delta-l_1)t}y(t)=
y(0)e^{\int_{-\infty}^{0}(H(s)-l_1)ds}.
$$
To prove (\ref{eq3.7}), it suffices to show that $\alpha<-\infty$.
From (\ref{eq3.3}) and (\ref{eq3.5}), a direct computation shows
that
$$H'(s)=-\frac{(a+1-b)p}{(p-1)(N-(a+1-b)p)}H(s)^{2-p}\xi(H(s)),$$
where
$$
\xi(s)=(p-1)s^p-(N-(a+1)p)s^{p-1}+\mu.
$$
From the definitions of $l_1,l_2$, we may suppose that
$$
H'(s)=(H(s)-l_1)(H(s)-l_2)g(H(s)),
$$
where $g$ is a continuous negative function on the interval
$[l_1,\delta]$, thus satisfies $|g(H(s))|\geqslant c_1>0$. From
(\ref{eq3.10}), it follows that
$$
\alpha=\underset{t\rightarrow-\infty}{\lim}e^{(\delta-l_1)t}y(t)=y(0)e^{\int_{-\infty}^{0}(H(s)-l_1){\rm
d}s} =y(0)e^{\int_{l_1}^{\delta}[(H(s)-l_2)g(H(s))]^{-1}{\rm
d}H(s)}.
$$
Since $l_2>\delta$ and $|g(H(s))|\geqslant c_1$ on $[l_1,\delta]$,
we know that
$$
\int_{l_1}^{\delta}[(H(s)-l_2)g(H(s))]^{-1}{\rm d}H(s)<+\infty,
$$
that is, $\alpha<+\infty$, thus (\ref{eq3.7}) follows. \epf

In the following corollary, we rewrite these conclusions on $y$
into those on the positive solution $u\in
\mathfrak{D}^{1,p}(\mathbb{R}^N)$ of equation (\ref{eq3.1}).

\begin{corollary}\label{coro3.1} Let $u\in
\mathfrak{D}^{1,p}(\mathbb{R}^N)$ be a positive solution of
equation (\ref{eq3.1}). Then there exists two positive constants
$C_1, C_2>0$ such that
\begin{equation}\label{eq3.11}
\lim_{r\rightarrow 0}\, r^{l_1}u(r)=C_1>0, \quad\quad
\lim_{r\rightarrow+\infty}\, r^{l_2}u(r)=C_2>0.
\end{equation}
and
\begin{equation}\label{eq3.13}
\lim_{r\rightarrow0}\, r^{l_1+1}|u'(r)|=C_1 l_1>0,
\quad\quad \lim_{r\rightarrow+\infty}\,
r^{l_2+1}|u'(r)|=C_2 l_2>0.
\end{equation}
 \end{corollary}
\proof From (\ref{eq3.2}), we know $u(r)=r^{-\delta}y(t)$.
Applying Lemma \ref{lem3.4} directly, we have
$$
\lim_{r\rightarrow0}\,
r^{l_1}u(r)=\lim_{t\rightarrow-\infty}\,
e^{(l_1-\delta)t}y(t)=y(0)c_1=C_1>0,
$$
$$
\lim_{r\rightarrow+\infty}\, r^{l_2}
u(r)=\lim_{t\rightarrow+\infty}\,
e^{(l_2-\delta)t}y(t)=y(0)c_2=C_2>0.
$$

Noting that $\lim\limits_{t\rightarrow-\infty}H(t)=l_1$ and
$\lim\limits_{t\rightarrow+\infty}\,H(t)=l_2$, it follows that
\begin{equation}\label{eq3.14}
\begin{split}
\lim_{r\rightarrow0}\,r^{l_1}u(r)\cdot H(t)
&=\lim_{r\rightarrow0}\,r^{l_1}u(r)\cdot
\frac{|z(t)|^{\frac{1}{p-1}}}{y(t)}
=\lim_{r\rightarrow0}\,r^{l_1}u(r)\cdot
\frac{r^{1+\delta}|u'(r)|}{r^{\delta}u(r)}\\
&=\lim_{r\rightarrow0}\,r^{l_1+1}|u'(r)| =C_1l_1>0
\end{split}
\end{equation}
and
\begin{equation}\label{eq3.15}
\begin{split}\lim_{r\rightarrow+\infty}\,r^{l_2}u(r)\cdot H(t)&
=\lim_{r\rightarrow+\infty}\,r^{l_2}u(r)\cdot
\frac{|z(t)|^{\frac{1}{p-1}}}{y(t)}
=\lim_{r\rightarrow+\infty}\,r^{l_2}u(r)\cdot
\frac{r^{1+\delta}|u'(r)|}{r^{\delta}u(r)}\\
&=\lim_{r\rightarrow+\infty}\,r^{l_2+1}|u'(r)|=C_2l_2>0.\
\end{split}
\end{equation}
\epf

Next, we shall give a uniqueness result of positive solution of
equation (\ref{eq3.1}).

\begin{theorem}\label{thm3.1} Suppose that $u_1(r)$ and $u_2(r)$
are two positive solutions of equation (\ref{eq3.1}). Let
$(y_1(t),z_1(t))$ and $(y_2(t),z_2(t))$ be two solutions to ODE
system (\ref{eq3.5}) corresponding to $u_1(r)$ and $u_2(r)$
respectively. If
\begin{equation}\label{eq3.16}
\underset{t\in(\infty,+\infty)}{\max}y_1(t)=y_1(0)=[\frac{N}{N-(a+1-b)p}(\delta^p-\mu)]^{\frac{1}{p_\ast-p}},
\end{equation}
and $y_2(0)=y_1(0)$. Then $(y_1(t),z_1(t))=(y_2(t),z_2(t))$, hence
$u_1=u_2$.
 \end{theorem}

\proof  The proof is similar to that of Theorem 3.11 in
\cite{AFP}. \epf

Similar to Theorem 3.13 in \cite{AFP}, we resume the above results
together and obtain the following theorem which describes the
asymptotic behavior of all the radial solutions to equation
(\ref{eq3.1}).

\begin{theorem}\label{thm3.2} All positive radial solutions to equation
(\ref{eq2.1}) have the form:
\begin{equation}\label{eq3.17}
u(\cdot)=\varepsilon^{-\frac{N-(a+1)p}{p}}u_0(\frac{\cdot}{\varepsilon}),
\end{equation}
where $u_0$ is a solution to equation (\ref{eq2.1}) satisfying
$u_0(1)=y(0)=[\frac{N}{N-(a+1-b)p}(\delta^p-\mu)]^{\frac{1}{p_\ast-p}}$.
Furthermore, there exist constants $C_1,C_2>0$ such that
\begin{equation}\label{eq3.18}
0<C_1\leqslant\frac{u_0(x)}{(|x|^{l_1/\delta}+|x|^{l_2/\delta})^{-\delta}}\leqslant
C_2,
\end{equation}
where $l_1,l_2$ are the two zeros of function
$\xi(s)=(p-1)s^p-(N-(a+1)p)s^{p-1}+\mu$ satisfying $0<l_1<l_2.$
 \end{theorem}

%%%%%%%%%%%%%%%%%%%%%%%%%%%%%%%%%%%%%%%%%%%%%%%%%%%%%%%%%%%%%%%%%%%%%%%%%%%%%%%%%%%%%%%%%%%%%%%%%%%%%

\section{Strict inequality (\ref{eq1.4})}
%%%%%%%%%%%%%%%%%%%%%%%%%%%%%%%%%%%%%%%%%%%%%%%%%%%%%%%%%%%%%%%%%%%%%%%%%%%%%%%%%%%%%%%%%%%%%%%%%%%%%

In this section, applying the asymptotic behavior of the solutions
to equation (\ref{eq2.1}) obtained in the previous section, we
give some estimates on the extremal function of $S_{0,\,\mu}$. Let
$u_0$ be an extremal function of $S_{0,\,\mu}$ with $\|u_0;L_b^{p_\ast}(\mathbb{R}^N)\|=1$. From the discussion
in Section 2 and 3, we know that $u_0$ is radial, and for all $\varepsilon>0$,
$$U_\varepsilon(r)=\varepsilon^{-\frac{N-(a+1)p}{p}}u_0(\frac{r}{\varepsilon})$$
is also an  extremal function of $S_{0,\,\mu}$, and there exists a
positive constant $C_\varepsilon$ such that $C_\varepsilon
U_\varepsilon$ is a solution to equation (\ref{eq2.1}). In fact,
from the proof of Lemma \ref{lem2.3}, we know that
$C_\varepsilon=S_{0,\,\mu}^{\frac{1}{p_\ast-p}}$, which is
independent of $\varepsilon$, denoted by $C_0$.  Set
$u_{\varepsilon}^{\ast}=C_0 U_\varepsilon$, then from equation
(\ref{eq2.1}) we have
\begin{equation}\label{eq4.1}
Q_\mu(u_{\varepsilon}^\ast)=\|u_{\varepsilon}^{\ast};
{L_b^{p_\ast}} \|^{p_\ast}=S_{0,\,\mu}^
{\frac{p_\ast}{p_\ast-p}}=S_{0,\,\mu}^{\frac{N}{(a+1-b)p}}.
\end{equation}

For any $\varepsilon>0$, and $m\in\mathbb{N}$ large enough such
that $B_{\frac{1}{m}}\subseteq \Omega$, define
\begin{equation}\label{eq4.2}
u_{\varepsilon}^m(x)= \Big\{{\begin{array}{l}
u_{\varepsilon}^\ast(x)-u_{\varepsilon}^\ast(\frac{1}{m}),\quad x\in B_{\frac{1}{m}}\backslash\{0\};\\
0,\quad\quad\quad\quad\quad\quad\quad\ \  x\in\Omega\backslash
B_{\frac{1}{m}}.
\end{array}}
\end{equation}

\begin{lemma}\label{lem4.1}
Set $\varepsilon=m^{-h},\ h>1$. Then as $m\to \infty$, we have
\begin{equation}\label{eq4.3}
Q_\mu(u_{\varepsilon}^m)\leqslant
S_{0,\,\mu}^{\frac{N}{(a+1-b)p}}+{\cal O}(m^{-(h-1)[(a+1+l_2)p-N]}),
\end{equation}
and
\begin{equation}\label{eq4.4}
\|u_{\varepsilon}^{\ast}; {L_b^{p_\ast}} \|^{p_\ast}\geqslant
S_{0,\,\mu}^{\frac{N}{(a+1-b)p}}- {\cal
O}(m^{-(h-1)[(b+l_2)p_*-N]}),
\end{equation}
where and afterward ${\cal O}(m^{-\alpha})$ denotes a positive quality which is
$O(m^{-\alpha})$, but is not $o(m^{-\alpha})$, as $m\to \infty$.
 \end{lemma}
\proof  We shall only prove (\ref{eq4.3}), and omit the prove of
(\ref{eq4.4}).

Since $Q_\mu(u_{\varepsilon}^m)=
\int_{\mathbb{R}^N}\frac{|Du_{\varepsilon}^m|^p}{|x|^{ap}}\,{\rm
d}x-\mu\int_{\mathbb{R}^N}\frac{|u_{\varepsilon}^m|^p}{|x|^{(a+1)p}}\,{\rm
d}x$, we estimate each term in $Q_\mu(u_{\varepsilon}^m)$ as
follows:
\begin{equation}\label{eq4.5}\begin{array}{rl}
\displaystyle \int_\Omega\frac{|{\rm
D}u_{\varepsilon}^m|^p}{|x|^{ap}}\,{\rm d}x &=\displaystyle
\int_{B_{\frac{1}{m}}}\frac{|{\rm
D}u_{\varepsilon}^\ast|^p}{|x|^{ap}}\,{\rm
d}x\\[5mm]
&=\displaystyle \int_{\mathbb{R}^N}\frac{|{\rm
D}u_{\varepsilon}^\ast|^p}{|x|^{ap}}\,{\rm d}x
-\int_{\mathbb{R}^N\backslash B_{\frac{1}{m}}}\frac{|{\rm
D}u_{\varepsilon}^\ast|^p}{|x|^{ap}}\,{\rm d}x\\
&\leqslant \displaystyle \int_{\mathbb{R}^N}\frac{|{\rm
D}u_{\varepsilon}^\ast|^p}{|x|^{ap}}\,{\rm d}x
\end{array}\end{equation}
and
\begin{equation}\label{eq4.6}\begin{array}{rl}
\displaystyle
&\displaystyle\int_\Omega\frac{|u_{\varepsilon}^m|^p}{|x|^{(a+1)p}}\,{\rm
d}x =\displaystyle \int_{B_{\frac{1}{m}}}
\frac{(u_{\varepsilon}^\ast(x)- u_{\varepsilon}^\ast(\frac1m)
)^p}{|x|^{(a+1)p}} \,{\rm d}x\\[5mm]
&\  \geqslant\displaystyle \int_{B_{\frac{1}{m}}}
\frac{u_{\varepsilon}^\ast(x)^p  -p
u_{\varepsilon}^\ast(\frac1m)u_{\varepsilon}^\ast(x)^{p-1}}{|x|^{(a+1)p}}
\,{\rm
d}x\\[5mm]
&\   =\displaystyle \int_{\mathbb{R}^N}
\frac{u_{\varepsilon}^\ast(x)^p}{|x|^{(a+1)p}} \,{\rm d}x
-\int_{\mathbb{R}^N\backslash
B_{\frac{1}{m}}}\frac{u_{\varepsilon}^\ast(x)^p}{|x|^{(a+1)p}}
\,{\rm d}x
-p u_{\varepsilon}^\ast(\frac1m) \displaystyle
\int_{B_{\frac{1}{m}}}\frac{u_{\varepsilon}^\ast(x)^{p-1}}{|x|^{(a+1)p}}
\,{\rm d}x.
\end{array}\end{equation}
On the other hand, from the definition of $u_{\varepsilon}^\ast$,
we have
\begin{equation}\label{eq4.7}\begin{array}{rl}
\displaystyle \int_{\mathbb{R}^N\backslash
B_{\frac{1}{m}}}\frac{u_{\varepsilon}^\ast(x)^p}{|x|^{(a+1)p}}
\,{\rm d}x &=C_0^p\omega_N \displaystyle \int_{\frac1m}^{+\infty}
\frac{\varepsilon ^{-[N-(a+1)p]} u_0(\frac
r\varepsilon)^p}{r^{(a+1)p}} r^{N-1}\,{\rm d}r\\[3mm]
& =C_0^p\omega_N \displaystyle \int_{m^{h-1}}^{+\infty}  u_0(t)^p
t^{N-1-(a+1)p} \,{\rm d}t
\\[3mm]
& ={\cal O}(m^{-(h-1)[(a+1+l_2)p-N]}),
\end{array}\end{equation}
where in the second equality, we make the change of variable
$t=\frac r\varepsilon$, and in the last equality, we use the
asymptotic behavior of $u_0$ at the infinity, since $h>1$, hence
$m^{h-1}\to \infty$ as $m\to \infty$. Note that
$\xi^\prime(l_2)=p(p-1)l_2^{p-1}-(p-1)(N-(a+1)p)l_2^{p-2}>0$, that
is $(a+1+l_2)p-N>0$. Similarly, we can estimate the last
integration in (\ref{eq4.6}) as follows:
\begin{equation}\label{eq4.8}\begin{array}{rl}
u_{\varepsilon}^\ast(\frac1m) \displaystyle
\int_{B_{\frac{1}{m}}}\frac{u_{\varepsilon}^\ast(x)^{p-1}}{|x|^{(a+1)p}}
\,{\rm d}x &=C_0^p\omega_N  u_0(\frac1{m\varepsilon})
\displaystyle \int_0^{\frac1m} \frac{\varepsilon ^{-[N-(a+1)p]}
u_0(\frac r\varepsilon)^{p-1}}{r^{(a+1)p}} r^{N-1}\,{\rm d}r\\[3mm]
&=C_0^p\omega_N  u_0(m^{h-1}) \displaystyle \int_0^{m^{h-1}}
u_0(t)^{p-1} t^{N-1-(a+1)p} \,{\rm d}t
\\[3mm]
& \leqslant C_0^p\omega_N  C_2 m^{-(h-1)l_2p}
[C+m^{(h-1)[N-(a+1)p-(p-1)l_2]} ]
\\[3mm]
&={\cal O}(m^{-(h-1)[(a+1+l_2)p-N]}),
\end{array}\end{equation}
where the last equality is from $\xi(l_2)=0$ and so
$N-(a+1)p-(p-1)l_2=\mu/l_2^{p-1}>0$. Thus, (\ref{eq4.3}) follows
from (\ref{eq4.5})-(\ref{eq4.8}). \epf

\begin{lemma}\label{lem4.2}
Set $\varepsilon=m^{-h},\ h>1$. If $c<(a+1+l_2)p-N$, then
\begin{equation}\label{eq4.9}
\int_{\mathbb{R}^N}\frac{|u_{\varepsilon}^m(x)|^p}{|x|^{(a+1)p-c}}\,{\rm
d}x\geqslant {\cal O}(m^{-ch}).
\end{equation}
 \end{lemma}

\proof A direct computation shows that
$$
\begin{array}{rl}
& \displaystyle
\int_{\mathbb{R}^N}\frac{|u_{\varepsilon}^m(x)|^p}{|x|^{(a+1)p-c}}\,{\rm
d}x =\displaystyle \int_{B_{\frac{1}{m}}}
\frac{(u_{\varepsilon}^\ast(x)- u_{\varepsilon}^\ast(\frac1m)
)^p}{|x|^{(a+1)p-c}} \,{\rm d}x\\[5mm]
&\ \ \  \geqslant\displaystyle \int_{B_{\frac{1}{m}}}
\frac{u_{\varepsilon}^\ast(x)^p  -p
u_{\varepsilon}^\ast(\frac1m)u_{\varepsilon}^\ast(x)^{p-1}}{|x|^{(a+1)p-c}}
\,{\rm
d}x\\[5mm]
& \ \ \  =\displaystyle \int_{\mathbb{R}^N}
\frac{u_{\varepsilon}^\ast(x)^p}{|x|^{(a+1)p-c}} \,{\rm d}x
-\int_{\mathbb{R}^N\backslash
B_{\frac{1}{m}}}\frac{u_{\varepsilon}^\ast(x)^p}{|x|^{(a+1)p-c}}
\,{\rm d}x -p u_{\varepsilon}^\ast(\frac1m) \displaystyle
\int_{B_{\frac{1}{m}}}\frac{u_{\varepsilon}^\ast(x)^{p-1}}{|x|^{(a+1)p-c}}
\,{\rm d}x.
\end{array}$$
We estimate each of the above integrations as follows:
\begin{equation}\label{eq4.10}
\displaystyle \int_{\mathbb{R}^N}
\frac{u_{\varepsilon}^\ast(x)^p}{|x|^{(a+1)p-c}} \,{\rm d}x
=C_0^p\omega_N\varepsilon^c \int_0^\infty
u_0(t)^pt^{N-1-(a+1)p+c}\, {\rm d}x
 = {\cal O}(m^{-ch}),
\end{equation}
\begin{equation}\label{eq4.11}
\begin{array}{rl}
\displaystyle \int_{\mathbb{R}^N\backslash
B_{\frac{1}{m}}}\frac{u_{\varepsilon}^\ast(x)^p}{|x|^{(a+1)p-c}}
\,{\rm d}x&=C_0^p\omega_N
\varepsilon^c\displaystyle\int_{m^{h-1}}^\infty
u_0(t)^pt^{N-1-(a+1)p+c}\, {\rm d}x\\
&={\cal O}(m^{-(h-1)[(a+1+l_2)p-N]-c})
\end{array}\end{equation}
and
\begin{equation}\label{eq4.12}
\begin{array}{rl}
u_{\varepsilon}^\ast(\frac1m) \displaystyle
\int_{B_{\frac{1}{m}}}\frac{u_{\varepsilon}^\ast(x)^{p-1}}{|x|^{(a+1)p-c}}
\,{\rm d}x &=C_0^p\omega_N  u_0(\frac1{m\varepsilon})
\displaystyle \int_0^{\frac1m} \frac{\varepsilon ^{-[N-(a+1)p]}
u_0(\frac r\varepsilon)^{p-1}}{r^{(a+1)p-c}} r^{N-1}\,{\rm d}r\\[3mm]
& =C_0^p\omega_N  u_0(m^{h-1}) \varepsilon^c\displaystyle
\int_0^{m^{h-1}} u_0(t)^{p-1} t^{N-1-(a+1)p+c} \,{\rm d}t
\\[3mm]
&\leqslant C_0^p\omega_N  C_2 m^{-(h-1)l_2p-ch}
[C+m^{(h-1)[N-(a+1)p-(p-1)l_2]} ]
\\[3mm]
& ={\cal O}(m^{-(h-1)[(a+1+l_2)p-N]-c}).
\end{array}\end{equation}
Note that since $c<(a+1+l_2)p-N$, we have $-ch>
-(h-1)[(a+1+l_2)p-N]-c$, that is, we prove the lemma.
 \epf

Let $\Omega$ be a smooth bounded open domain in $\mathbb{R}^N$
with $0\in \Omega$, define $\mathfrak{D}_{a,b}^{1,p}(\Omega)$ as
the closure of $C_0^\infty (\Omega)$ under the norm $\| u
\|_{\mathfrak{D}_{a,b}^{1,p}(\Omega)}=\| |{\rm D}u|; L_a^p
(\Omega) \|$ and
\begin{equation}\label{eq4.13}
S_{\lambda,\,\mu}(p,a,b;\Omega)=\inf\left\{Q_{\lambda,\,\mu}(u)\ :\
u\in \mathfrak{D}_{a,b}^{1,p}(\Omega),\
 \|u;L_b^{p_\ast}(\Omega)\|=1\right\},
\end{equation}
where
$$
Q_{\lambda,\,\mu}(u)=\int_{\Omega}\frac{|{\rm
D}u|^p}{|x|^{ap}}\,{\rm
d}x-\mu\int_{\Omega}\frac{|u|^p}{|x|^{(a+1)p}}\,{\rm d}x - \lambda
\int_{\Omega}\frac{|u|^p}{|x|^{(a+1)p-c}}\,{\rm d}x.
$$
If $\lambda=0$, by rescaling argument, it is easy to show that
$S_{0,\,\mu}(p,a,b;\Omega)=S_{0,\,\mu}$. But for $\lambda>0$, we shall
have a strict inequality between $S_{\lambda,\,\mu}(p,a,b;\Omega)$
and $S_{0,\,\mu}$.

\begin{theorem}\label{thm4.1}
If $\mu\in (0,\ \overline{\mu}),\ \lambda>0,\ b\in [a,\ a+1),\ c\in (0,\
(a+1+l_2)p-N)$, then the strict inequality (\ref{eq1.4}) holds.
\end{theorem}

\proof We shall study
$$
\frac{Q_{\lambda,\,\mu}(u_\varepsilon^m)}{\|u_\varepsilon^m;
L_b^{p_*}(\Omega)\|^p}.
$$
It follows from Lemma 4.1 and 4.2 that
\begin{equation}\label{eq4.14}
\begin{array}{rl}
Q_{\lambda,\,\mu}(u_\varepsilon^m) &=Q_\mu(u_\varepsilon^m)
-\lambda\displaystyle \int_\Omega
\frac{|u_{\varepsilon}^m(x)|^p}{|x|^{(a+1)p-c}}\,{\rm
d}x \\
&\leqslant S_{0,\,\mu}^{\frac{N}{(a+1-b)p}}+{\cal
O}(m^{-(h-1)[(a+1+l_2)p-N]})-{\cal O}(m^{-ch})
\end{array}\end{equation}
and
\begin{equation}\label{eq4.15}
\begin{array}{rl}
\|u_\varepsilon^m; L_b^{p_*}(\Omega)\|^p &\geqslant
  S_{0,\,\mu}^{\frac{N}{(a+1-b)p_*}} -{\cal
O}(m^{-(h-1)[(b+l_2)p_*-N]p/p_*})\\
  &= S_{0,\,\mu}^{\frac{N}{(a+1-b)p_*}}-{\cal
O}(m^{-(h-1)[(a+1+l_2)p-N]}).
\end{array}\end{equation}
Thus, we have
\begin{equation}\label{eq4.16}
\begin{array}{rl}
\dfrac{Q_{\lambda,\,\mu}(u_\varepsilon^m)}{\|u_\varepsilon^m;
L_b^{p_*}(\Omega)\|^p} &\leqslant
 \dfrac{S_{0,\,\mu}^{\frac{N}{(a+1-b)p}}+{\cal
O}(m^{-(h-1)[(a+1+l_2)p-N]})-{\cal
O}(m^{-ch})}{S_{0,\,\mu}^{\frac{N}{(a+1-b)p_*}} -{\cal
O}(m^{-(h-1)[(b+l_2)p_*-N]p/p_*})}\\
  &= S_{0,\,\mu}  +{\cal
O}(m^{-(h-1)[(a+1+l_2)p-N]}) - {\cal O}(m^{-ch}).
\end{array}\end{equation}
If $c\in (0,\ (a+1+l_2)p-N)$, we can choose $h$ large enough such
that $c<(h-1)(a+1+l_2)p-N)/h$ and so $-ch > -(h-1)[(a+1+l_2)p-N]$,
thus as $m$ large enough, (\ref{eq1.4}) holds. \epf

%%%%%%%%%%%%%%%%%%%%%%%%%%%%%%%%%%%%%%%%%%%%%%%%%%%%%%%%%%%%%%%%%%%%%%%%%%%%%%%%%%%%%%%%%%
\section{Application}
%%%%%%%%%%%%%%%%%%%%%%%%%%%%%%%%%%%%%%%%%%%%%%%%%%%%%%%%%%%%%%%%%%%%%%%%%%%%%%%%%%%%%%%%%%
In this section, as an application of the strict inequality of (\ref{eq1.4}),
we consider the existence of nontrivial solutions to the following quasilinear Brezis-Nirenberg type problem involving Hardy potential
 and Sobolev critical exponent:

\begin{equation}\label{eq5.1}
\left\{ {\begin{array}{rl}  -\mbox{div}(\dfrac{|{\rm D}u|^{p-2}{\rm
D}u}{|x|^{ap}})-\mu\dfrac{|u|^{p-2}u}{|x|^{(a+1)p}}
&=\dfrac{|u|^{{p_\ast}-2}u}{|x|^{bp_\ast}}+\lambda\dfrac{|u|^{p-2}u}{|x|^{(a+1)p-c}},   \mbox{ in}\ \Omega,\\[3mm]
u&=0,  \mbox{ on}\ \partial\,\Omega,
\end{array}} \right.
\end{equation}
where $\Omega\subset \mathbb{R}^N$ is an open bounded domain with
$C^1$ boundary and $0\in \Omega$, $1<p<N,\
p_\ast=\frac{Np}{N-(a+1-b)p}$, $0\leq a<\frac{N-p}{p},\ a\leqslant b<(a+1),\
c>0$; $\lambda,\ \mu$ are two positive real parameters.

To obtain the existence result, let's define the energy functional $E_{\lambda,\,\mu}$ on
$\mathfrak{D}_{a,b}^{1,p}(\Omega)$ as
$$
E_{\lambda,\,\mu}(u)=\frac1p\int_\Omega \left[  \frac{|{\rm D}u|^p}{|x|^{ap}}-\mu \frac{|u|^p}{|x|^{(a+1)p}}
-\lambda\frac{|u|^p}{|x|^{(a+1)p-c}}\right]\,{\rm
d}x -\frac1 {p_*}\int_\Omega \frac{|u|^{p_*}}{|x|^{bp_*}}\,{\rm
d}x.
$$
It is easy to see that $E_{\lambda,\,\mu}$ is well-defined in $\mathfrak{D}_{a,b}^{1,p}(\Omega)$,
and $E_{\lambda,\,\mu} \in C^1(\mathfrak{D}_{a,b}^{1,p}(\Omega), \mathbb{R})$. Furthermore, all the critical points of $E_{\lambda,\,\mu}$
 are weak solutions to (\ref{eq5.1}). We shall apply the Mountain Pass Lemma without (PS) condition due to Ambrosetti and Rabinowitz
 \cite{AR} to ensure the existence of (PS)$_\beta$ sequence of $E_{\lambda,\,\mu}$ at some
 Mountain Pass type minimax value level $\beta$. Then the strict inequality (\ref{eq1.4})
 implies that $\beta< \frac {a+1-b}N S_{0, \mu}^\frac{N}{(a+1-b)p}$. Finally, combining the generalized concentration compactness
 principle and a compactness property called singular Palais-Smale condition due to Boccardo and Murat \cite{BM}(cf. also \cite{GP}),
 we shall obtain the existence of nontrivial solutions to (\ref{eq5.1}).

Let's define two more functionals on $\mathfrak{D}_{a,b}^{1,p}(\Omega)$ as follows:
$$
I_\mu(u)=\frac1p\int_\Omega \frac{|{\rm D}u|^p}{|x|^{ap}}\,\,{\rm
d}x -\frac\mu p\int_\Omega \frac{|u|^p}{|x|^{(a+1)p}}\,\,{\rm
d}x,\ J(u)=\int_\Omega \frac{|u|^p}{|x|^{(a+1)p-c}}\,\,{\rm
d}x,
$$
and denote ${\cal M}=\{u\in \mathfrak{D}_{a,b}^{1,p}(\Omega)\ : \ J(u)=1 \}$. For $\mu\in (0, \overline{\mu})$, the Hardy inequality shows that
$\frac1p \frac{|{\rm D}u|^p}{|x|^{ap}}\,{\rm
d}x -\frac\mu p \frac{|u|^p}{|x|^{(a+1)p}}\,{\rm
d}x$ is nonnegative measure on $\Omega$. The classical results in the Calculus of Variations(cf. \cite{SM})
show that $I_\mu$ is lower semicontinuity on ${\cal M}$.
On the other hand the compact imbedding theorem in \cite{XBJ} implies that  ${\cal M}$ is weakly closed.
Thus the direct method ensure that $I_\mu$ attains its minimum on ${\cal M}$,
denote $\lambda_1=\min \{I_\mu(u) \ :\ u\in {\cal M} \}>0$. From the homogeneity of $I_\mu$ and $J$, $\lambda_1$ is the first
nonlinear eigenvalue of problem:
\begin{equation}\label{eq5.3}
\left\{ {\begin{array}{rl} -\mbox{div}(\dfrac{|{\rm D}u|^{p-2}{\rm
D}u}{|x|^{ap}})-\mu\dfrac{|u|^{p-2}u}{|x|^{(a+1)p}}
&=\lambda\dfrac{|u|^{p-2}u}{|x|^{(a+1)p-c}},   \mbox{in}\ \Omega,\\
u&=0,  \mbox{on}\ \partial\,\Omega.
\end{array}} \right.
\end{equation}

The following lemma indicates that $E_{\lambda,\,\mu}$ satisfies the geometric condition of Mountain Pass Lemma without (PS) condition due to Ambrosetti and Rabinowitz
 \cite{AR}, the proof is direct and omitted.

 \begin{lemma}\label{lem5.1}
If $\mu\in (0,\ \mu), \lambda\in (0,\ \lambda_1)$, then
\begin{enumerate}
\item[(i)]$E_{\lambda,\,\mu}(0)=0$;

\item[(ii)] $\exists \,\alpha, r>0$, s.t. $E_{\lambda,\,\mu}(u)\geqslant \alpha$, if $\|u\|=r$;

\item[(iii)] For any $v\in \mathfrak{D}_{a,b}^{1,p}(\Omega),\ v\neq 0$, there exists $T>0$ such that
$E_{\lambda,\,\mu}(tv)\leqslant 0$ if $t>T$.
\end{enumerate}
\end{lemma}

For $v\in \mathfrak{D}_{a,b}^{1,p}(\Omega)$ with $\|v\|>r$ and $E_{\lambda,\,\mu}(v)\leqslant 0$,
set
$$\beta:=\inf_{\gamma\in\Gamma} \max_{t\in [0, 1]} E_{\lambda,\,\mu}(\gamma(t)),$$
where
$$
\Gamma:=\{\gamma\in C([0, 1], \mathfrak{D}_{a,b}^{1,p}(\Omega)) \ | \ \gamma(0)=0,\ \gamma(1)=v \}.
$$
It is easy to see that $\beta$ is independent of the choice of $v$ such that $E_{\lambda,\,\mu}(v)\leqslant 0$,
and furthermore $\beta\geqslant \alpha$. If $\beta$ is finite, from Lemma \ref{lem5.1} and Mountain Pass Lemma, there exists a (PS)$_\beta$ sequence $\{u_m\}_{m=1}^\infty$ of $E_{\lambda,\,\mu}$ at level $\beta$,
that is, $E_{\lambda,\,\mu} (u_m) \to \beta$ and $E_{\lambda,\,\mu}^\prime (u_m) \to 0$ in the dual space $(\mathfrak{D}_{a,b}^{1,p}(\Omega))^\prime$
of $\mathfrak{D}_{a,b}^{1,p}(\Omega)$ as $m\to \infty$.

 \begin{lemma}\label{lem5.2}
If $\mu\in (0,\ \overline{\mu}), \lambda\in (0,\ \lambda_1)$, then the strict inequality (\ref{eq1.4}) is equivalent to
\begin{equation}\label{eq5.4}\beta< \frac{a+1-b}{N}S_{0,\,\mu}^{\frac{N}{(a+1-b)p}}.\end{equation}
\end{lemma}
\proof \textbf{1.} (\ref{eq1.4}) $\implies$ (\ref{eq5.4}).

Let $v_1$ be a function such that $\|v_1; L_b^{p_*}(\Omega)\|=1$, and $Q_{\lambda,\,\mu}(v_1)<S_{0,\,\mu}$. We have
\begin{equation}\label{eq5.5}
 \begin{array}{rl}\beta & \leqslant\sup\limits _{0<t<\infty} E_{\lambda,\,\mu}(tv_1)=\sup\limits  _{0<t<\infty}
 (\dfrac{t^p}p Q_{\lambda,\,\mu}(v_1)-\dfrac{t^{p_*}}{p_*})\\
& = (\dfrac 1p-\dfrac1{p_*})Q_{\lambda,\,\mu}(v_1)^{\frac{p_*}{p_*-p}}=\dfrac{a+1-b}N Q_{\lambda,\,\mu}(v_1)^{\frac{N}{(a+1-b)p}}\\
& <\dfrac{a+1-b}N S_{0,\,\mu}^{\frac{N}{(a+1-b)p}}.
\end{array}
\end{equation}

\textbf{2.} (\ref{eq5.4}) $\implies$ (\ref{eq1.4}).

Since $\lambda<\lambda_1$, for $u=g(t)=tv$ with $t$ closed to $0$, we have $(D E_{\lambda,\,\mu}(u), u)>0$;
while for $u=g(1)=v$, we have
$$
(D E_{\lambda,\,\mu}(v), v)< p E_{\lambda,\,\mu}(v)\leqslant 0.
$$
Consider function $f(t)=E_{\lambda,\,\mu}(tv)\in C^1([0,1],\, \mathbb{R})$, we have that $f^\prime (t)>0$ for $t$ closed to $0$,
and $f_-^\prime(1)\leqslant 0$. From the medium value theorem, there exists $t_0\in (0, 1)$ such that $f^\prime(t_0)=0$, that is,
for $u=t_0v$, we have
$$
(D E_{\lambda,\,\mu}(u), u)=Q_{\lambda,\,\mu}(u)-\|u; L_b^{p_*}(\Omega)\|^{p_*}=0.
$$
Thus a direct computation shows that
$$
\frac{Q_{\lambda,\,\mu}(u)}{\|u; L_b^{p_*}(\Omega)\|^{p}} =Q_{\lambda,\,\mu}(u)^{1-p/p_*}
=(\frac{N}{a+1-b} E_{\lambda,\,\mu}(u) )^\frac{(a+1-b)p}{N},
$$
that is,
$$
\beta=\inf_{\gamma\in\Gamma} \max_{t\in [0, 1]} E_{\lambda,\,\mu}(\gamma(t))\geqslant
\frac{a+1-b}N S_{\lambda,\,\mu}(p, a, b, \Omega)^\frac{N}{(a+1-b)p}.
$$
Hence (\ref{eq5.4}) $\implies$ (\ref{eq1.4}). \epf

\begin{lemma}\label{lem5.3}
If $\mu\in (0,\ \overline{\mu}), \lambda\in (0,\ \lambda_1)$, then any (PS)$_\beta$ sequence of $E_{\lambda,\,\mu}$ is bounded in
$\mathfrak{D}_{a,b}^{1,p}(\Omega)$.
\end{lemma}
\proof Suppose that $\{u_m\}_{m=1}^\infty$ is a (PS)$_\beta$ sequence of $E_{\lambda,\,\mu}$. As $m\to \infty$, we have
\begin{equation}\label{eq5.6}
 \begin{array}{ll} & \beta+o(1) =E_{\lambda,\,\mu}(u_m)\\[2mm]
 &\ \ \ \ =\dfrac1p\displaystyle \int_\Omega \left[  \frac{|{\rm D}u_m|^p}{|x|^{ap}}-\mu \frac{|u_m|^p}{|x|^{(a+1)p}}
-\lambda\frac{|u_m|^p}{|x|^{(a+1)p-c}}\right]\,{\rm
d}x -\dfrac1 {p_*}\displaystyle\int_\Omega \frac{|u_m|^{p_*}}{|x|^{bp_*}}\,{\rm
d}x
\end{array}\end{equation}
and
\begin{equation}\label{eq5.7}
 \begin{array}{rl} & o(1)\|\varphi\| =({\rm D} E_{\lambda,\,\mu}(u_m), \varphi)\\[2mm]
 &\ \ \ \ \ =\displaystyle\int_\Omega \left[  \frac{|{\rm D}u_m|^{p-2} {\rm D}u_m \cdot {\rm D}\varphi}{|x|^{ap}}
 -\mu \frac{|u_m|^{p-2} u_m \varphi}{|x|^{(a+1)p}}
-\lambda\frac{|u_m|^{p-2} u_m \varphi}{|x|^{(a+1)p-c}}\right]\,{\rm
d}x \\[4mm]
&\ \ \  \ \ \ \ \ \ -\displaystyle\int_\Omega \frac{|u_m|^{p_*-2}u_m \varphi}{|x|^{bp_*}}\,{\rm
d}x,
\end{array}\end{equation}
for any $\varphi \in \mathfrak{D}_{a,b}^{1,p}(\Omega)$. From (\ref{eq5.6}) and (\ref{eq5.7}), as $m\to \infty$, it follows that
$$
\begin{array}{rl} p_*\beta+o(1)-o(1)\|u_m\|& =p_*E_{\lambda,\,\mu}(u_m) - ({\rm D} E_{\lambda,\,\mu}(u_m), u_m)\\[2mm]
 &=(\dfrac{p_*}p-1) \displaystyle \int_\Omega \left[  \frac{|{\rm D}u_m|^p}{|x|^{ap}}-\mu \frac{|u_m|^p}{|x|^{(a+1)p}}
-\lambda\frac{|u_m|^p}{|x|^{(a+1)p-c}}\right]\,{\rm
d}x  \\[4mm]
& \geqslant (\dfrac{p_*}p-1)(1-\dfrac{\lambda}{\lambda_1})\displaystyle \int_\Omega \left[
\frac{|{\rm D}u_m|^p}{|x|^{ap}}-\mu \frac{|u_m|^p}{|x|^{(a+1)p}}
\right]\,{\rm
d}x\\[4mm]
& \geqslant (\dfrac{p_*}p-1)(1-\dfrac{\lambda}{\lambda_1})(1-\dfrac{\mu}{\overline{\mu}})\|u_m\|^p.
\end{array}
$$
Thus, $\{u_m\}_{m=1}^\infty$ is bounded in
$\mathfrak{D}_{a,b}^{1,p}(\Omega)$ if $\mu\in (0,\ \mu), \lambda\in (0,\ \lambda_1)$. \epf

From the boundedness of $\{u_m\}_{m=1}^\infty$ in
$\mathfrak{D}_{a,b}^{1,p}(\Omega)$, we have the following medium convergence:
$$
u_m \rightharpoonup u  \mbox{ in }   \mathfrak{D}_{a,b}^{1,p}(\Omega), \ L_1^{p}(\Omega)  \mbox{ and } L_b^{p_*}(\Omega),
$$
$$
u_m \to u  \mbox{ in }   L_\alpha^{r}(\Omega)\ \mbox{ if } 1\leq r<\frac{Np}{N-p},\ \frac\alpha r< (a+1)+N(\frac1r-\frac1p),
$$
$$
u_m \to u\ \  \mbox{a.e. in } \Omega.
$$

In order to obtain the strong convergence of $\{u_m\}_{m=1}^\infty$ in $L_b^{p_*}(\Omega)$, we need the following generalized concentration compactness
principle(cf. also \cite{TY}) and \cite{WW} and references therein), the proof is similar to that in \cite{LPL2} and we omit it.

\begin{lemma}[Concentration Compactness Principle]
\label{lem5.4}
Suppose that ${\cal M}(\mathbb{R}^N)$ is the space of bounded measures on $\mathbb{R}^N$, and
$\{u_m\}\subset \mathfrak{D}_{a,b}^{1,p}(\Omega)$ is a sequence such that:
$$
\begin{array}{ll}
u_m \rightharpoonup u & \mbox{ in }   \mathfrak{D}_{a,b}^{1,p}(\Omega),\\[2mm]
\xi_m:=\left(|x|^{-ap}|{\rm D}u_m|^p-\mu  |x|^{-(a+1)p}|u_m|^p \right) \,{\rm
d}x \rightharpoonup \xi & \mbox{ in }   {\cal M}(\mathbb{R}^N),\\[2mm]
\nu_m:=|x|^{-bp_*}|u_m|^{p_*}  \,{\rm
d}x\rightharpoonup \nu & \mbox{ in }   {\cal M}(\mathbb{R}^N),\\[2mm]
u_m\to u & \mbox{ a.e. on }   \mathbb{R}^N.
\end{array}
$$
Then there are the following statements:
\begin{enumerate}
\item[(1)] There exists some at most countable set $J$, a family $\{x^{(j)}\ :\ j\in J\}$ of distinct
points in $\mathbb{R}^N$, and a family $\{\nu^{(j)}\ :\ j\in J \}$ of positive numbers such that
\begin{equation}
\label{eq5.8}
\nu=|x|^{-bp_*}|u|^{p_*}  \,{\rm
d}x+\sum_{j\in J} \nu^{(j)}\delta_{x^{(j)}},
\end{equation}
where $\delta_x$ is the Dirac-mass of mass $1$ concentrated at $x\in \mathbb{R}^N$.

\item[(2)] The following inequality holds
\begin{equation}
\label{eq5.9}
\xi \geq (|x|^{-ap}|{\rm D}u|^p-\mu  |x|^{-(a+1)p}|u|^p ) \,{\rm
d}x+\sum_{j\in J} \xi^{(j)}\delta_{x^{(j)}},
\end{equation}
for some family $\{\xi^{(j)}>0\ :\ j\in J \}$ satisfying
\begin{equation}
\label{eq5.10}
S_{0,\, \mu}\big( \nu^{(j)}\big)^{p/p_*}\leqslant \xi^{(j)},\ \ \mbox{for all }j\in J.
\end{equation}
In particular, $\sum\limits_{j\in J}\big( \nu^{(j)}\big)^{p/p_*}<\infty$.
\end{enumerate}
\end{lemma}

\begin{lemma}\label{lem5.5} If $\mu\in (0,\ \mu), \lambda\in (0,\ \lambda_1)$, let $\{u_m\}_{m=1}^\infty$ be a (PS)$_\beta$ sequence of $E_{\lambda,\,\mu}$ at level $\beta$
defined above. (\ref{eq5.4}) implies that $\nu^{(j)}=0$ for all $j\in J$, that is, up to a subsequence,
 $u_m\to u$ in $L_b^{p_*}(\Omega)$ as $m\to 0$.
\end{lemma}

\proof From Lemma \ref{lem5.3}, $\{u_m\}_{m=1}^\infty$ is bounded in $\mathfrak{D}_{a,b}^{1,p}(\Omega)$, then we have that
$|{\rm D}u_m|^{p-2} {\rm D}u_m$ is bounded in $\left(L^{p\prime}(\Omega; |x|^{-ap}) \right)^N$,
where $p^\prime$ is the conjugate exponent of $p$, i.e. $\frac1p+\frac1{p^\prime}=1$.
Without loss of generality, we suppose that
$T\in \left(L^{p\prime}(\Omega; |x|^{-ap}) \right)^N$ such that
$$
|{\rm D}u_m|^{p-2} {\rm D}u_m \rightharpoonup  T \mbox{  in  }\left(L^{p\prime}(\Omega; |x|^{-ap}) \right)^N.
$$

Also, $|u_m|^{p-2} u_m$ is bounded in $L^{p\prime}(\Omega; |x|^{-(a+1)p})$,
$|u_m|^{p_*-2} u_m$ is bounded in $L^{p_*\prime}(\Omega; |x|^{-bp_*})$, and $u_m\to u$ almost everywhere in $\Omega$,
thus it follows that
$$
|u_m|^{p-2} u_m \rightharpoonup  |u|^{p-2} u  \mbox{  in  } L^{p\prime}(\Omega; |x|^{-(a+1)p})
$$
and
$$
|u_m|^{p_*-2} u_m \rightharpoonup  |u|^{p_*-2} u  \mbox{  in  } L^{p_*\prime}(\Omega; |x|^{-bp_*}).
$$
From the compactness imbedding theorem in \cite{XBJ}, it follows that
$$
|u_m|^{p-2} u_m \to  |u|^{p-2} u  \mbox{  in  } L^{p\prime}(\Omega; |x|^{-(a+1)p+c}).
$$
Taking $m\to \infty$ in (\ref{eq5.7}), we have
\begin{equation}
\label{eq5.14}
\displaystyle\int_\Omega  \frac{T\cdot {\rm D}\varphi}{|x|^{ap}}\,{\rm
d}x=
\mu \displaystyle\int_\Omega \frac{|u|^{p-2} u \varphi}{|x|^{(a+1)p}}\,{\rm
d}x
+\lambda \displaystyle\int_\Omega \frac{|u|^{p-2} u \varphi}{|x|^{(a+1)p-c}}\,{\rm
d}x + \displaystyle\int_\Omega \frac{|u|^{p_*-2}u \varphi}{|x|^{bp_*}}\,{\rm
d}x,
\end{equation}
for any $\varphi \in \mathfrak{D}_{a,b}^{1,p}(\Omega)$. Let $\varphi=\psi u_m$ in (\ref{eq5.7}),
where $\psi \in C(\bar\Omega)$, and take $m\to \infty$, it follows that
\begin{equation}
\label{eq5.15}
\displaystyle\int_\Omega \psi \,{\rm
d}\xi+ \displaystyle\int_\Omega \frac{uT\cdot {\rm D}\psi}{|x|^{ap}}\,{\rm
d}x=
\displaystyle\int_\Omega \psi \,{\rm
d}\nu+
 \lambda \displaystyle\int_\Omega \frac{|u|^{p}  \psi}{|x|^{(a+1)p-c}}\,{\rm
d}x.
\end{equation}
Let $\varphi=\psi u$ in (\ref{eq5.14}), it follows that
\begin{equation}
\label{eq5.16}\begin{array}{ll}
\displaystyle\int_\Omega  \frac{uT\cdot {\rm D}\psi}{|x|^{ap}}\,{\rm
d}x& + \displaystyle\int_\Omega  \frac{\psi T\cdot {\rm D}u}{|x|^{ap}}\,{\rm
d}x =
\mu \displaystyle\int_\Omega \frac{|u|^{p}  \psi}{|x|^{(a+1)p}}\,{\rm
d}x\\[3mm]
 & \ \ \
+\lambda \displaystyle\int_\Omega \frac{|u|^{p} \psi}{|x|^{(a+1)p-c}}\,{\rm
d}x+ \displaystyle\int_\Omega \frac{|u|^{p_*} \psi}{|x|^{bp_*}}\,{\rm
d}x,
\end{array}\end{equation}
Thus, form (\ref{eq5.8}) and Lemma \ref{lem5.4}, (\ref{eq5.15})$-$(\ref{eq5.16}) implies that
\begin{equation}
\label{eq5.17}\begin{array}{ll}
\displaystyle\int_\Omega \psi \,{\rm
d}\xi &=\displaystyle\int_\Omega  \frac{\psi T\cdot {\rm D}u}{|x|^{ap}}\,{\rm
d}x-\mu \displaystyle\int_\Omega \frac{|u|^{p}  \psi}{|x|^{(a+1)p}}\,{\rm
d}x+\displaystyle\int_\Omega \psi \,{\rm
d}\nu-\displaystyle\int_\Omega \frac{|u|^{p_*} \psi}{|x|^{bp_*}}\,{\rm
d}x\\[3mm]
& =\displaystyle\int_\Omega  \frac{\psi T\cdot {\rm D}u}{|x|^{ap}}\,{\rm
d}x- \mu \displaystyle\int_\Omega \frac{|u|^{p}  \psi}{|x|^{(a+1)p}}\,{\rm
d}x+\sum_{j\in J} \nu^{(j)}\psi(x^{(j)}).
\end{array}\end{equation}
Letting $\psi\to \delta_{x^{(j)}}$, we have
$$
\xi^{(j)}= \nu^{(j)}.
$$
Combining with (\ref{eq5.10}), it follows that $\nu^{(j)} \geqslant S_{0,\, \mu}\big( \nu^{(j)}\big)^{p/q}$, which means that
\begin{equation}
\label{eq5.18}\nu^{(j)} \geqslant S_{0,\, \mu}^\frac{N}{(a+1-b)p},
\end{equation}
 if $\nu^{(j)}\neq 0$. On the other hand, taking $m\to\infty$ in (\ref{eq5.6}), and using (\ref{eq5.17}) with $\psi\equiv 1$,
(\ref{eq5.8}) and (\ref{eq5.14}),  it follows that
\begin{equation}
\label{eq5.19}\begin{array}{ll}
 \beta & = \dfrac1p \displaystyle\int_\Omega \,{\rm
d}\xi -\dfrac{1}{p_*}\displaystyle\int_\Omega \,{\rm
d}\nu - \dfrac\lambda p \displaystyle\int_\Omega \frac{|u|^{p} }{|x|^{(a+1)p-c}}\,{\rm
d}x  \\[3mm]
& =\dfrac1p\left( \sum\limits_{j\in J} \nu^{(j)}
+\displaystyle\int_\Omega  \frac{T\cdot {\rm D}u}{|x|^{ap}}\,{\rm
d}x- \mu \displaystyle\int_\Omega \frac{|u|^{p}  }{|x|^{(a+1)p}}\,{\rm
d}x\right)\\[3mm]
&\ \  \ -\dfrac1{p_*}\left( \sum\limits_{j\in J} \nu^{(j)}
 + \displaystyle\int_\Omega \frac{|u|^{p_*} }{|x|^{bp_*}}\,{\rm
d}x\right)- \dfrac\lambda p \displaystyle\int_\Omega \frac{|u|^{p} }{|x|^{(a+1)p-c}}\,{\rm
d}x  \\[3mm]
&=(\dfrac1p- \dfrac1{p_*}) \sum\limits_{j\in J} \nu^{(j)}
+ (\dfrac1p- \dfrac1{p_*}) \displaystyle\int_\Omega \frac{|u|^{p_*}}{|x|^{bp_*}}\,{\rm
d}x\\[3mm]
& \geqslant (\dfrac1p- \dfrac1{p_*}) \sum\limits_{j\in J} \nu^{(j)}
=\dfrac {a+1-b}N \sum\limits_{j\in J} \nu^{(j)}.
\end{array}\end{equation}
From (\ref{eq5.18}), (\ref{eq5.19}), (\ref{eq5.4}) implies that $\nu^{(j)}=0$ for all $j\in J$. Hence we have
$$
\int_\Omega \frac{|u_m|^{p_*} }{|x|^{bp_*}}\,{\rm
d}x \to \int_\Omega \frac{|u|^{p_*} }{|x|^{bp_*}}\,{\rm
d}x,
$$
as $m\to \infty$. Thus, the Brezis-Lieb Lemma \cite{BL} implies that, up to a subsequence,
 $u_m\to u$ in $L_b^{p_*}(\Omega)$ as $m\to 0$. \epf

In order to deduce the almost everywhere convergence of ${\rm D}u_m$ in $\Omega$ and to obtain existence of nontrivial solution
to (\ref{eq5.1}), we shall apply the variational approach supposed in \cite{GP}
and a convergence theorem due to Boccardo and Murat(cf. Theorem 2.1 in \cite{BM}), so we suppose that $a=0$,
and $\mathfrak{D}_{a,b}^{1,p}(\Omega) =W_0^{1,p}(\Omega)$.

\begin{theorem}\label{thm5.1}
If $a=0, \mu\in (0,\ \overline{\mu}),\ \lambda\in (0,\ \lambda_1),\ b\in [0,\ 1),\ c\in (0,\
(1+l_2)p-N)$, then there exists a nontrivial solution to (\ref{eq5.1}).
\end{theorem}
\proof Apply the variational approach supposed in \cite{GP} and a convergence theorem in \cite{BM}, there exists a subsequence
of $\{u_m\}_{m=1}^\infty$, still denoted by $\{u_m\}_{m=1}^\infty$, such that
$$
u_m\to u \mbox{  in  } W_0^{1,\,q}(\Omega),\ q<p,
$$
which implies that $u$ is a solution to (\ref{eq5.1}) in sense of distributions. Since $u\in W_0^{1,p}(\Omega)$,
by density argument, $u$ is a weak solution to (\ref{eq5.1}). Next, we shall show that $u\not\equiv 0$.

In fact, from the homogeneity and Lemma \ref{lem5.5}, we have
$$
\begin{array}{ll}
0<\alpha\leqslant  \beta & =\lim\limits_{m\to\infty}E_{\lambda,\,\mu}(u_m)
=\lim\limits_{m\to\infty}\left[E_{\lambda,\,\mu}(u_m)
-\dfrac1p({\rm D} E_{\lambda,\,\mu}(u_m), u_m)\right]\\
&=\lim\limits_{m\to\infty}(\dfrac1p- \dfrac1{p_*}) \displaystyle\int_\Omega \frac{|u_m|^{p_*}}{|x|^{bp_*}}\,{\rm
d}x\\[3mm]
&=(\dfrac1p- \dfrac1{p_*}) \displaystyle\int_\Omega \frac{|u|^{p_*}}{|x|^{bp_*}}\,{\rm
d}x,
\end{array}
$$
Thus, $u\not\equiv 0$. \epf

In sight of Theorem \ref{thm5.1}, we conjecture that the conclusion is also true for $0\leqslant a <\frac{N-p}{p}$.
\begin{conjecutre}
If $0\leqslant a <\frac{N-p}{p}, \mu\in (0,\ \overline{\mu}),\ \lambda\in (0,\ \lambda_1),\ b\in [a,\ a+1),\ c\in (0,\
(a+1+l_2)p-N)$, then there exists a nontrivial solution to (\ref{eq5.1}).
\end{conjecutre}

%%%%%%%%%%%%%%%%%%%%%%%%%%%%%%%%%%%%%%%%%%%%%%%%%%%%%%%%%%%%%%%%%%%%%%%%%%%%%%%%%%%%%%%%%%


\begin{thebibliography}{01}

\bibitem{AFP}B.Abdellaoui, V.Felli and I.Peral, Existence and Nonexistence Results for Quasilinear Elliptic Equations
Involving the p-laplacian, Advances in Differential Equations,
481-508.
\bibitem{AR}A. Ambrosetti \& P. H. Rabinowitz, Dual variational methods in critical point theory and applications,
J. Funct. Anal., \textbf{14}(1973), pp349-381.
\bibitem{BM}L. Boccardo and F. Murat, Almost everywhere convergence of the gradients of solutions to elliptic and parabolic equations,
Nonli. Anal., TMA, Vol. {\bf 19}(1992), 581-597.
\bibitem{BL}H. Brezis and E. Lieb, A relation between pointwise convergence of functions and convergence of functionals, Proc. Amer. Math. Soc., Vol. {\bf 88}(1983), PP486-490.
\bibitem{BN}H. Brezis and L. Nirenberg, Positive solutions of nonlinear elliptic equations involving critical exponents,
Comm. Pure Appl. Math., Vol. {\bf 36}(1983), 437-477.
\bibitem{CKN}L. Caffarrelli, R. Kohn and L. Nirenberg, First order interpolation inequalities with weights,
Compositio Mathematica, Vol. {\bf 53}(1984), 259-275.
\bibitem{CH}D. M. Cao and P. G. Han, Solutions for semilinear elliptic equations with critical exponents and Hardy potential,
J. Diff. Eqns., Vol. {\bf 205}(2004), 521-537.
\bibitem{CG}K.-S. Chou and D. Geng, On the critical dimension of a semilinear degenerate elliptic equation involving critical Sobolev-Hardy exponent,
Nonli. Anal., TMA, Vol. {\bf 26}(1996), PP1965-1984.
\bibitem{EH1}H. Egnell, Semilinear elliptic equations involving critical Sobolev exponents, Arch. Rational Mech. Anal., Vol. {\bf 104}(1988), PP27-56.
\bibitem{EH2}H. Egnell, Existence and nonexistence results for m-Laplace equations involving critical Sobolev exponents, Arch. Rational Mech. Anal., Vol. {\bf 104}(1988), PP57-77.

\bibitem{FG}A. Ferrero and F. Gazzola, Existence of solutions for singular critial growth semilinear elliptic equations,
J. Diff. Eqns., Vol. {\bf 177}(2001), 494-522.

\bibitem{GP}J. P. Garcia Azorero and I. Peral Alonso, Hardy inequalities and some critical elliptic and parabolic problems,
J. Diff. Eqns., Vol. {\bf 144}(1998), 441-476.

\bibitem{GV}M. Guedda and L. Veron, Quasilinear ellptic equations involving critical Sobolev exponents, Nonli. Anal., TMS, Vol. {\bf 13}(1989), PP879-902.

\bibitem{HT}T. Horiuchi, Best constant in weighted Sobolev inequality with weights being powers of distance from the origen, J. Inequal. Appl., Vol. {\bf 1}(1997), PP275-292.

\bibitem{JE}E. Jannelli, The role played by space dimension in elliptic critical problems,
J. Diff. Eqns., Vol. {\bf 156}(1999), 407-426.


\bibitem{JS}E. Jannelli and S. Solomini, Critical behaviour of some elliptic equations with singular potentials, Rapport no. 41/96, Dipartimento di Mathematica Universita degi Studi di Bari, 70125 Bari, Italia.


\bibitem{LPL2}P. L. Lions, The concentration-compactness principle in the calculus of variations, the limit case,
Rev. Mat. Ibero Americana, part 1, Vol. {\bf 1}(1985), PP145-201, part 2, Vol. {\bf 2}(1985), PP45-121.


\bibitem{NL}L. Nicolaescu, A weighted semilinear elliptic equation involving critical Sobolev exponents, Diff. \& Int. Eqns., Vol. {\bf 3}(1991), PP653-671.

\bibitem{PS}P. Pucci and J. Serrin, Critical exponents and critical dimensions for polyharmonic operators, J. Math. Pures Appl., Vol. {\bf 69}(1990), PP55-83.

\bibitem{RW}D. Ruiz and M. Willem, Elliptic problems with critical exponents and Hardy potentials,
J. Diff. Eqns., Vol. {\bf 190}(2003), PP524-538.


\bibitem{SM}M. Struwe, Variational Methods, Applications to Nonlinear Partial Differential Equations and Hamiltonian Systems, 2 ed, Springer-Verlag, 1996.


\bibitem{TY}J.-G. Tan and J.-F. Yang, On the singular variational problems,
Acta Mat. Sinica, Vol. {\bf 24}(2004), PP672-690.

\bibitem{WW}Z.-Q. Wang and M. Willem, Singular minimization problems,
J. Diff. Eqns., Vol. {\bf 161}(2000), PP307-320.


\bibitem{XC}B.-J. Xuan and Z.-C. Chen, Existence, multiplicity and bifurcation for critical polyharmonic equations, Sys. Sci. and Math. Sci., Vol. {\bf 12}(1999), PP59-69.


\bibitem{XBJ}B.-J. Xuan, The solvability of quasilinear
Brezis-Nirenberg type problems with singular weights, Nonli. Anal., Vol. {\bf 62}(2005), PP703-725.



\bibitem{XSY}B.-J. Xuan, S.-W. Su and Y.-J. Yan, Existence results of
Brezis-Nirenberg problems with Hardy potential and singular
coefficients, accepted by Nonlinear Analysis.

\bibitem{ZXP}X.-P. Zhu, Nontrivial solution of quasilinear elliptic equations involving critical Sobolev exponent, Sci. Sinica, Ser. A Vol. {\bf 31}(1988), PP1166-1181.




\end{thebibliography}
\end{document}